\documentstyle{amsppt}
\magnification 1200
\NoBlackBoxes

\pageheight{9 true in}
\pagewidth{6.5 true in}
\def\vphi{\varphi}
\def\hchi{\hat \chi}
\def\hsigma{\hat \sigma}

\topmatter
\title Upper bounds for $|L(1,\chi)|$
\endtitle
\author
Andrew Granville and K. Soundararajan
\endauthor

\address
Department of Mathematics, University of Georgia, Athens, GA , USA
\endaddress
\email
andrew\@math.uga.edu
\endemail
\address
School of Mathematics,  Institute for Advanced Study,
Princeton, NJ 08540, USA
\endaddress
\email
ksound\@math.ias.edu
\endemail
\thanks{The first author is supported, in part, by the National Science
Foundation.  The second author
is partially supported by the American Institute
of Mathematics (AIM)}
\endthanks

\endtopmatter

\head 1. Introduction \endhead

\noindent Given a non-principal Dirichlet character $\chi \pmod q$,
an important problem
in number theory is to obtain good estimates for the size of $L(1,\chi)$.
The best bounds known give that $q^{-\epsilon} \ll_{\epsilon} |L(1,\chi)|
\ll \log q$, while assuming the Generalized Riemann Hypothesis,
J.E. Littlewood showed that $1/\log\log q \ll |L(1,\chi)|\ll \log\log q$.
Littlewood's result reflects the true range of the size of $|L(1,\chi)|$
as it is known that there exist characters $\chi_\pm$
for which $L(1,\chi_+)\asymp \log\log q$ and $L(1,\chi_-)\asymp 1/\log\log q$.

In this paper we focus on sharpening the upper bounds known for
$|L(1,\chi)|$; in particular, we wish to determine constants $c$ (as
small as possible) for which
the bound $|L(1,\chi)| \le (c+o(1)) \log q$ holds.
To set this in context,
observe that if $X$ is such that $\sum_{n\le x} \chi(n) = o(x)$ for all
$x>X$ then
$$
L(1,\chi) = \sum_{n\le X} \frac{\chi(n)}{n} + o(\log q). \tag{1.1}
$$
Trivially $X=q$ is permissible and so $|L(1,\chi)| \le (1+o(1))\log q$.
Less trivially the P{\' o}lya-Vinogradov inequality gives that
$X=q^{\frac 12 +o(1)}$ is permissible.  Finally note that
D. Burgess' character sums estimates permit one to take
$X=q^{\frac 14 +o(1)}$ if $q$ is cube-free, and $X=q^{\frac 13+o(1)}$
otherwise.  In particular we get that $|L(1,\chi)| \le (1/4+o(1))\log p$
when $q=p$ is prime.
In [1] Burgess improved on this ``trivial'' bound, for quadratic characters,
obtaining that $L(1,\fracwithdelims(){\cdot}{p})
\le 0.2456 \log p$ for all large primes $p$.  This was subsequently improved
by P.J. Stephens [7] to
$L(1,\fracwithdelims(){\cdot}{p}) \le (2-2/\sqrt{e} +o(1))
\frac{1}{4}\log p$.  This
result is best-possible in the sense that one can construct  totally  multiplicative functions $f$, taking only values $-1$ and $1$, such that
$\sum_{n\le x} f(n) = o(x)$ for all $x>X$, and
$\sum_{n\le X} f(n)/n \sim (2-2/\sqrt{e})\log X$.
Stephens' result was extended by Pintz [8] to all quadratic characters.
No analogous improvements over the trivial bound were known
for complex characters $\chi$.  We give such a result below.

\proclaim{Corollary}  Define  $c_2= 2-2/\sqrt e=0.786938\ldots,\
c_3 = 4/3 -1/e^{2/3} =0.819916\ldots, \ c_4 =0.8296539741\ldots$,
and $c_k = c_\infty := 34/35$ for $k\ge 5$.
For any primitive Dirichlet character $\chi \pmod q$ of order $k$, we have
$$
|L(1,\chi)| \le
\cases
\frac{1}{4} (c_k +o(1)) \log q \qquad&\text{if  } q \text{  is cube-free},\\
\frac 13 (c_k +o(1)) \log q \qquad &\text{otherwise}.\\
\endcases
$$
\endproclaim

We establish the Corollary by focussing more generally on
multiplicative functions satisfying a ``Burgess-type'' condition.
Given a subset $S$ of the unit disc ${\Bbb U}$, we define
${\Cal F}(S)$ to be the class of all completely multiplicative
functions $f$ such that $f(p)\in S$ for all primes $p$.
We denote by ${S}_k$ the set $\{ 0\} \cup \{ \xi:\ \xi^k=1\}$.
Our problem is to bound $(1/\log X) |\sum_{n\le X} f(n)/n|$
for $f\in {\Cal F}(S)$ assuming that
$$
\sum_{n\le x} f(n) = o(x), \tag{1.2}
$$
for suitable $x\ge X$.
More precisely, let $A \ge 1$ be a parameter, and define
$$
\gamma(S; A) :=
\limsup_{X\to \infty} \max\Sb f\in
{\Cal F}(S)\\ \text{(1.2) holds for } X\le x\le X^A\endSb
\frac{1}{\log X} \Big| \sum_{n\le X} \frac{f(n)}{n}\Big|.
$$
Note that for fixed $S$, $\gamma(S;A)$ is a non-increasing function of
$A$.  Further if $S_1 \supset S_2$ then $\gamma(S_1;A) \ge \gamma(S_2;A)$
for all $A\ge 1$.  Set $\gamma(S) =\lim_{A\to \infty} \gamma(S;A)$.

By (1.1) and Burgess' estimates we see that
 if $\chi$ is a character $\pmod q$ of order $k$ then
$$
|L(1,\chi)| \le
\cases
\frac{1}{4}
(\gamma(S_k) +o(1)) \log q \qquad&\text{if  } q \text{  is cube-free},\\
\frac 13 (\gamma(S_k) +o(1)) \log q \qquad &\text{otherwise}.\\
\endcases
$$
Thus our corollary above follows from our main Theorem which
establishes upper bounds on $\gamma(S;A)$.

\proclaim {Theorem 1} With the definitions as above, $\gamma(S_k;1) \le
c_k$ for $k=2$, $3$ and $4$.  For $k \ge 5$ we
have
$$
\gamma(S_k) \le \gamma(S_k;\sqrt{e}) \le \gamma({\Bbb U};\sqrt{e})
\le c_{\infty}.
$$
\endproclaim

It is possible to show that $\gamma({\Bbb U};1) \le c$ for an absolute
constant $c<1$.  This follows from
P.D.T.A. Elliott's groundbreaking result [2]
that the magnitude of averages of multiplicative functions
varies slowly.  Precisely, for
any $f\in {\Cal F}({\Bbb U})$ and $1\le w\le x$, we have
the following Lipschitz-type estimate
$$
\frac 1x \Big| \sum_{n\le x} f(n)\Big| -\frac wx \Big| \sum_{n\le x/w}
f(n)\Big| \ll \left( \frac{\log 2w}{\log x}\right)^{\frac 1{19}}.
$$
Thus if $\sum_{n\le X} f(n) =o(X)$, then there is some $\delta<1$ such
that for all $X^\delta \le x\le X$  we have $|\sum_{n\le x} f(n)| \le x/2$.
Hence
$$
\align
\frac{1}{\log X} \Big| \sum_{n\le X} \frac{f(n)}{n} \Big|
&= \frac{1}{\log X}
\Big| \int_{1}^{X} \sum_{n\le t} f(n) \frac{dt}{t^2} \Big|
+O\left(\frac{1}{\log X}\right)\\
&\le \frac 1{\log X} \left( \int_1^{X^\delta} \frac{dt}{t} +
\int_{X^\delta}^{X}
\frac{dt}{2t} \right) +o(1) = c +o(1).
\\
\endalign
$$
where $c=(1+\delta)/2$.
Elliott's exponent $1/19$ has recently been improved in [4] to any
exponent $< 1-2/\pi$,
and is probably true for any exponent $<1$.
However it seems that the value of $c$
given by this method is inevitably much closer to $1$ than $c_\infty$.

Although we do not go into this here, one can, via Lipschitz-type estimates,
improve the upper bound from Theorem 1 for $k=3$ to
$\gamma(S_3;1)<c_3-\delta$ for some tiny $\delta>0$.

By means of a construction we are also able to
give lower bounds for $\gamma(S)$.

\proclaim{Theorem 2a} We have
$\gamma({\Bbb U}) \ge \gamma(S_{2k}) \ge \gamma(S_2) \ge (2-2/\sqrt{e})$.
If $k$ is odd then $\gamma(S_k) \ge (1+\delta_k) (1-e^{-1/(1+\delta_k)})$
where $\delta_k= \cos (\pi/k)$.
\endproclaim

Combining with Theorem 1 we get that $\gamma(S_2) = 2-2/\sqrt{e}$.
It is tempting to conjecture that $\gamma({\Bbb U}) = 2-2/\sqrt{e}$.

Returning to our application to bounding $L(1,\chi)$, we note that
we have not exploited all the information on characters
available to us.  Namely, if $\chi$ is a character of order $k$ then
$\chi^j$ is a non-principal character for $j=1$, $2$, $\ldots k-1$,
so that the Burgess estimates apply to mean values of $\chi^j$ as well.
Although we have not been able to take advantage of this fact, we
can establish some limits on how much it can imply.
The problem is to bound $(1/\log X) |\sum_{n\le X} f(n)/n|$ for
a given $f\in {\Cal F}(S_k)$ satisfying
$$
\sum_{n\le x} f(n)^j = o(x)  \ \ \text{\rm for } 1\leq j\leq k-1, \tag{1.3}
$$
for suitable $x \ge X$.
Precisely, for $A\ge 1$ we wish to determine
$$
\gamma_k (A) := \limsup_{X\to \infty} \max\Sb f\in  {\Cal F}(S_k)\\
\text{(1.3) holds for } X\le x\le X^A \endSb
\frac{1}{\log X} \Big| \sum_{n\le X} \frac{f(n)}{n}\Big|.
$$
Plainly $\gamma_k(A)$ is a decreasing function of $A$, and $\gamma_k(A)
\le \gamma(S_k;A)$.  We set $\gamma_k = \lim_{A\to \infty} \gamma_k(A)$.
If $\chi \pmod q$ is a character of order $k$ then
$$
|L(1,\chi)| \le
\cases
\frac{1}{4} (\gamma_k +o(1)) \log q \qquad&\text{if  } q \text{  is cube-free},\\
\frac 13 (\gamma_k +o(1)) \log q \qquad &\text{otherwise}.\\
\endcases
$$

\proclaim{Theorem 2b} For large $k$ we have
$$
\gamma_k \ge (e^{\gamma}+ o_k(1)) \frac{\log \log k}{\log k}.
$$
\endproclaim

We prove Theorem 2b in Section 6; indeed there we shall establish
a more precise lower bound on $\gamma_k$, and give numerical
data for small $k$.  We suspect that Theorem 2b gives the
correct size of $\gamma_k$ for large $k$.  At any rate, it
seems safe to conjecture that $\gamma_k = o_k(1)$, which
would imply that for any fixed $\epsilon >0$ and if $k$ is sufficiently large
then $|L(1,\chi)| \le \epsilon \log q$ for all characters $\chi\pmod q$ of
order $k$.

Our final result obtains an upper bound for $|L(1,\chi)|$ on ``average''
over all the characters of order $k$.

\proclaim{Theorem 3}  Suppose that $f\in S_k$ satisfies (1.4)
for all $X\le x\le X^{\vphi(k)+1}$.  Then
$$
\biggl \{ \prod\Sb 1\leq j\leq k-1\\ (j,k)=1\endSb
\frac{1}{\log X} \Big| \sum_{n\le X} \frac{f(n)^j}{n}\Big|
\biggr\}^{1/\vphi(k)} \leq \Big\{  \frac{43}{15} e^\gamma +o_k(1) \Big\}
\ \frac{\log\log k}{\log k}.
$$
Consequently
$$
\biggl\{ \prod\Sb \chi \pmod q \\ \chi \ \text{\rm of order} \ k\endSb
|  L(1,\chi) | \biggr\}^{1/\vphi(k)} \le
\cases
 \left\{  \frac{43}{60} e^\gamma +o_k(1) \right\}
\frac{\log\log k}{\log k}\ \log q \qquad&\text{if  } q \text{  is cube-free},\\
\left\{  \frac{43}{45} e^\gamma +o_k(1) \right\}
\frac{\log\log k}{\log k}\ \log q \qquad &\text{otherwise}.\\
\endcases
$$
\endproclaim

{\sl Acknowledgements.}  We are grateful to Roger Heath-Brown for
some stimulating conversations on this topic.

\head 2.  Preliminaries \endhead

\noindent Define $y:= \exp((\log X)^{\frac 14})$.  In proving Theorem
1  it is convenient to restrict attention to completely multiplicative
functions $f$ satisfying $f(p)=1$ for all $p\le y$.  We indicate
first why this entails no loss in generality.

\proclaim{Lemma 1}  Let $f$ be a multiplicative function
with $|f(n)|\le 1$ for all $n$.  Then
$$
\frac{1}{\log X} \Big| \sum_{n\le X} \frac{f(n)}{n} \Big|
\ll \exp\Big( -\frac 12 \sum_{p\le X} \frac{1-\text{\rm Re }f(p)}{p}\Big).
$$
\endproclaim

\demo{Proof}  See Proposition 8.1, and the comments  following it, in [3].
\enddemo

In proving Theorem 1  we may thus assume that
$$
\sum_{p\le X} \frac{1-\text{Re } f(p)}{p} \ll 1.
$$
Since $|1-f(p)|^2 \ll (1-\text{Re }f(p))$, we deduce by the
Cauchy-Schwarz inequality that
$$
\sum_{p\le X} \frac{|1-f(p)|}{p}
\le \Big( \sum_{p\le X} \frac{1}{p}\Big)^{\frac 12} \Big( \sum_{p\le X}
\frac{|1-f(p)|^2}{p} \Big)^{\frac 12} \ll \sqrt{\log \log X}.
\tag{2.1}
$$

\proclaim{Lemma 2}  Suppose $f$ is a multiplicative
function with $|f(n)|\le 1$ for all $n$, and that
$f$ satisfies {\rm (2.1)}.  Let $f_s(n)$
be the completely multiplicative function defined by
$f_s(p)=f(p)$ if $p>y$ and $f_s(p)=1$ for $p\le y$.
Define
$$
\Theta(f,y) = \prod_{p\le y} \left(1-\frac 1p\right) \left(1+\frac{f(p)}{p}
+\frac{f(p^2)}{p^2} + \ldots \right).
$$
Then for all $X^2 \ge x>X$,
$$
\frac 1x \sum_{n\le x} f(n) = \Theta(f,y)
\frac 1x\sum_{n\le x} f_s(n) + O( (\log x)^{-\frac 12}),
$$
and
$$
\frac 1{\log X} \sum_{n\le X} \frac{f(n)}{n} =
\Theta(f,y)
\frac{1}{\log X} \sum_{n\le X} \frac{f_s(n)}{n} + O( (\log X)^{-\frac 14}).
$$
\endproclaim

\demo{Proof}  The first assertion
follows from (2.1) and Proposition 4.5 of [3], while
the second assertion follows from Proposition 8.2 of [3].
\enddemo

Note that $|\Theta(f,y)|\le 1$ always.
If $|\Theta(f,y)| =o(1)$ then the bound in Theorem 1 is
immediate.   If $|\Theta(f,y)| \gg 1$
and $f$ meets the hypothesis of Theorem 1, then
$f_s$ meets the hypothesis of Theorem 1,
and it suffices to demonstrate the conclusion for $f_s$.
Thus Lemma 2 allows us to restrict attention
to completely multiplicative
functions $f$ with $f(p)=1$ for all primes $p\le y$, and
we suppose this henceforth.

\proclaim{Lemma 3} Let $f$ be a completely multiplicative
function with $f(p)=1$ for all $p\le y$, and $|f(p)|\leq 1$ otherwise,
and let $g$ be the completely multiplicative function
defined by $g(p)=|1+f(p)|-1$.  Put $G(u) = \sum_{n\le u} g(n)$.
Then for any $y\le u\le X$
$$
\align
\frac {1}{\log X} \Big| \sum_{n\le X} \frac{f(n)}{n} \Big|
&\le \frac{1}{\log X} \sum_{n\le X} \frac{g(n)}{n} +o(1)
\\
&\le \frac{1}{\log u} \int_{1}^{u} \frac{|G(t)|}{t^2} dt  +o(1).
\\
\endalign
$$
\endproclaim

\demo{Proof}  Note that
$$
\sum_{n\le X} \frac{f(n)}{n} = \frac{1}{X} \sum_{n\le X} \sum_{d|n} f(d)
+ O(1).
$$
Now $|\sum_{d|n} f(d)| \le \sum_{d|n} g(d)$ unless
$n$ is divisible by the square of some prime $p$ with $f(p) \neq 1$, so
that $p>y$.
The contribution of such $n$ is readily bounded by
$$
\frac 1X \sum_{p\ge y} \sum\Sb n\le X\\ p^2 |n \endSb d(n) \ll \log X \sum_{p> y}
\frac{1}{p^2} \ll \frac{\log X}{y} .
$$
It follows that
$$
\frac{1}{\log X} \Big| \sum_{n\le X} \frac{f(n)}{n} \Big| \le
 \frac{1}{\log X} \sum_{n\le X} \frac{g(n)}{n} + o(1)
= \frac{1}{\log X} \int_1^X \frac{G(t)}{t^2} dt +o(1).
$$

  From Lemma 2.1 of [4] we know that
$$
\frac{|G(u)|}{u} \le \frac{1}{\log u} \int_1^u \frac{|G(t)|}{t^2} dt +
O\left(\frac{1}{\log u}\right).
$$
  From this it follows that
$$
\left( \frac{1}{\log u}\int_1^u \frac{|G(t)|}{t^2} dt \right)^{\prime}
= -\frac{1}{u\log^2 u} \int_1^u \frac{|G(t)|}{t^2}dt + \frac{|G(u)|}
{u^2 \log u} \le O\left(\frac{1}{u\log^2 u}\right).
$$
We deduce that
$$
\align
\frac{1}{\log X} \int_1^X \frac{|G(t)|}{t^2} dt
&\le \frac{1}{\log u} \int_1^u \frac{|G(t)|}{t^2} dt
+ O\left(\int_u^X \frac{1}{t\log^2 t} dt \right)
\\
&= \frac{1}{\log u} \int_1^u \frac{|G(t)|}{t^2} dt
+ O\left(\frac{1}{\log u}\right),
\\
\endalign
$$
showing that $(1/\log u)\int_1^u |G(t)|/t^2 dt$ is essentially a non-increasing function, and the Lemma follows.
\enddemo

We record the value of an integral
that we will encounter several times. For $C>0$, we have
$$
\frac 1{\log u} \int_{u^{e^{-C}}}^u \left( C - \log \left( \frac {\log u}{\log t} \right) \right) \frac{dt}t = C-1 + e^{-C}. \tag{2.2}
$$

\proclaim{Lemma 4}  Let $f$ and $g$ be as in Lemma 3,
and put $I(u) =\sum_{p\le u} (1-g(p))/p$.  If $I(u) \le 1$ then
$$
\frac{1}{\log u} \int_1^u \frac{|G(t)| }{t^2} dt \le \frac{1}{2}
+ \frac{1}{2\log u} \int_1^u  (1-I(t))(1-I(u/t)) \frac{dt}{t} + o(1),
\tag{2.3}
$$
Useful bounds on the right hand side of {\rm (2.3)}  are
$$
1 - \frac{I(\sqrt{u})}{2}+o(1)
, \qquad
 1 - \frac{(1-I(\sqrt{u}))}{\log u}\int_1^u I(t) \frac{dt}{t}+o(1),
\qquad 1-\frac{1}{\log u} \int_{\sqrt{u}}^{u} I(t) \frac {dt}{t}+o(1),
$$
and
$$
3-I(u)-2e^{-I(u)/2} +o(1).
$$
\endproclaim

\demo{Proof} If $f(p)=1$ for every prime $p$ for which
$p^2$ divides  $n$ then, by induction on the number of
primes dividing $n$, we see that
$$
g(n) \ge 1-\sum_{p|n} (1-g(p)), \qquad \text{and} \qquad
g(n)\le 1-\sum_{p|n} (1-g(p)) + \sum_{pq|n} (1-g(p))(1-g(q)).
$$
It follows that
$$
\align
G(t) &\ge \sum_{n\le t} \Big (1-\sum_{p|n} (1-g(p))\Big)
+ O\Big( \sum_{p\ge y} \sum
\Sb n\le t\\ p^2 |n \endSb 1\Big) \\  &= t - t\sum_{p\le t} \frac{1-g(p)}{p}
 +o(t)= t(1-I(t))+o(t), \tag{2.4} \\
\endalign
$$
and, similarly, that
$$
G(t) \le t (1-I(t)) + \frac{t}{2} \sum_{p\le t} \sum_{q\le t/p}
\frac{1-g(p)}{p} \frac{1-g(q)}{q} +o(t).
$$
Thus if $I(u)\le 1+o(1)$ then $G(t) \ge o(t)$ for all $t\le u$, and
so
$$
\int_{1}^{u} \frac{|G(t)|}{t^2} dt
\le \int_1^u \Big(1-I(t) +\frac{1}{2} \sum_{p\le t} \sum_{q\le t/p}
\frac{1-g(p)}{p} \frac{1-g(q)}{q} \Big) \frac{dt}{t} + o(\log t).
$$
Since
$$
\int_1^u I(t) \frac{dt}{t} = \frac 12 \int_1^u (I(t)+I(u/t)) \frac{dt}{t},
$$
and
$$
\int_1^u I(t) I(u/t) \frac{dt}{t} = \int_1^u \sum_{p\le t} \sum_{q\le t/p}
\frac{1-g(p)}{p}\frac{1-g(q)}{q} \frac{dt}{t}
$$
we obtain the upper bound (2.3).  Since both $I(t)$ and $I(u/t)$ are
in $[0,1]$ and one of them is at least $I(\sqrt{u})$ we
immediately obtain our first alternative bound on the
RHS of (2.3).  Next
$$
\align
\text{RHS of (2.3)} &= 1-\frac{1}{\log u} \int_1^u I(t) \frac{dt}{t}
+ \frac{1}{2\log u}\int_1^u I(t)I(u/t) \frac{dt}{t} \\
&\le 1- \frac{1}{\log u} \int_1^u I(t) \frac{dt}{t} +
\frac{1}{\log u} \int_{\sqrt{u}}^{u} I(t) I(\sqrt{u}) \frac{dt}{t},\\
\endalign
$$
which proves our second alternative bound.  Further
$$
\align
\frac{1}{2\log u} \int_1^u (1-I(t))(1-I(u/t)) \frac{dt}{t}
&= \frac{1}{\log u}\int_{\sqrt{u}}^{u} (1-I(t))(1-I(u/t))\frac{dt}{t}  \\
&\le \frac{1}{\log u} \int_{\sqrt{u}}^{u} (1-I(t))\frac{dt}{t}, \\
\endalign
$$
which gives our third alternative bound.  Lastly note that
$1-I(t)\ge 0$ and also
$$
1-I(t) = 1-I(u) + \sum_{t\le p\le u} \frac{1-g(p)}{p} \le
1-I(u) + 2\log \left(\frac{\log u}{\log t}\right) +o(1),
$$
so that using our third alternative bound we get that
the RHS of (2.3) is
$$
\le \frac 12+ \frac{1}{\log u} \int_{\sqrt{u}}^{u^{e^{-I(u)/2}}} \frac{dt}{t}
+ \frac{1}{\log u}\int_{u^{e^{-I(u)/2}}}^{u} \left(1-I(u)
+ 2\log \left(\frac{\log u}
{\log t}\right)\right)\frac{dt}{t} + o(1),
$$
from which our final bound follows by (2.2).

\enddemo

\proclaim{Lemma 5}  Assume that $f\in {\Cal F}(\Bbb U)$ satisfies {\rm (1.2)}
for $X< x\le Y\le X^2$. Define $h(n) =\sum_{ab=n} f(a) \overline{f(b)}$.
If $Y \le X$ then $|\sum_{n\le Y} h(n)| \le
Y\log Y +O(Y)$.  If $X< Y\le X^2$ then
$$
\Big|\sum_{n\le Y} h(n)\Big| \le Y \log \frac{X^2}{Y} + o(Y\log Y).
$$
\endproclaim

\demo{Proof}  Since $|h(n)|\le d(n)$ we see that $|\sum_{n\le Y}
h(n)| \le \sum_{n\le Y} d(n) = Y\log Y +O(Y)$ which gives the
first assertion.  Suppose now that $X< Y\le X^2$ and write
$$
\sum_{n\le Y} h(n) = \sum_{ab\le Y} f(a) \overline{f(b)}
= \left(\sum\Sb ab\le Y\\ a\le Y/X\endSb + \sum\Sb ab\le Y \\ b\le Y/X \endSb
+ \sum\Sb ab\le Y\\ a, b \ge Y/X\endSb - \sum\Sb ab\le Y\\ a, b
\le Y/X\endSb \right) f(a) \overline{f(b)}.
$$
If $a\le Y/X$ then $X^2 \ge Y/a>X$, and so by (1.2) the first sum above is
$$
\sum_{a\le Y/X} \sum_{b\le Y/a} f(a) \overline{f(b)}
= \sum_{a\le Y/X} o(Y/a) = o(Y\log Y).
$$
Similarly the second sum above $\sum_{ab\le Y, \ b\le Y/X}$
is also $o(Y\log Y)$.  The third term is (in magnitude)
$$
\le \sum_{Y/X \le a \le X} \sum_{Y/X \le b\le Y/a} 1 = Y \log \frac{X^2}{Y} +
o(Y\log Y).
$$
Finally the last term is $O(Y^2/X^2) = o(Y\log Y)$ since $Y\le X^2$.
The Lemma follows.
\enddemo

\head 3.  Proof of Theorem 1 for $k=2, 3$ and $4$ \endhead

\noindent
Throughout this section we shall only assume that (1.2) holds for $x=X$.

\subhead 3a.  Proof of Theorem 1 for  $k=2$ \endsubhead

\noindent Suppose $k=2$ and $f\in {\Cal F}({\Cal S}_2)$.  Note here
that $g(p)=f(p)$.  From (2.4) it follows that
$\sum_{n\le X} f(n) = G(X) \ge X(1-I(X)) +o(X)$, so that by (1.2) we get
$I(X) \ge 1+o(1)$.  Let $u$ be the largest integer $\le X$ with
$I(u) \le 1$.  Plainly $u\ge y$ is large.
Note that $I(u)=1+o(1)$, and so by Lemma 3 and the
final bound in Lemma 4 it follows that
$$
\frac{1}{\log X} \Big| \sum_{n\le X} \frac{f(n)}{n} \Big|
\le 3-I(u)-2 e^{-I(u)/2} + o(1) = \left(2-\frac{2}{\sqrt{e}}\right) + o(1).
$$

\subhead 3b.  Proof of Theorem 1 for $k=3$ \endsubhead

\noindent Here note that $g(p)=1$ if $f(p)=1$, and $g(p)=0$ if
$f(p)\neq 1$.  Also note that
$$
\text{\rm Re } f(n) \ge 1- \sum_{p|n} (1-\text{Re } f(p))
\ge 1-\frac 32 \sum_{p|n} (1-g(p)) .
$$
Hence
$$
o(X) =\sum_{n\le X} \text{Re }f(n) \ge
X - \frac {3X}2 \sum_{p\le X}\frac{1-g(p))}p + o(X)
= X \left\{ 1 - \frac 32 I(X) +o(1)\right\} ,
$$
so that $I(X) \ge 2/3 +o(1)$.  Let $u$ be the largest integer below $X$
with $I(u)\le 2/3$.  Note that $u\ge y$ is
large and that $I(u)=2/3+o(1)$.

For $t\le u$ note that $0\le 1-I(t) \le 1$ and that
$$
\align
1-I(t) &= \frac 13 +I(u)-I(t) +o(1) =\frac{1}{3} + \sum_{t\le p\le u}
\frac{1-g(p)}{p} +o(1) \le \frac{1}{3} +\sum_{t\le p\le u} \frac 1p +o(1) \\
&= \frac 13+ \log \left(\frac{\log u}{\log t}\right) +o(1).
\\
\endalign
$$
Hence, using (2.2),
$$
\align
\int_{\sqrt{u}}^{u} (1-I(t))\frac{dt}{t}
&\le \int_{\sqrt{u}}^{u^{e^{-2/3}}} \frac{dt}{t} + \int_{u^{e^{-2/3}}}^{u}
\left(\frac 13+ \log \left( \frac{\log u}{\log t}\right) \right)\frac{dt}{t} +o(\log u) \\
&= \left(\frac 56 -\frac{1}{e^{2/3}} +o(1)   \right) \log u .\\
\endalign
$$
Applying Lemma 3 and the third bound in Lemma 4  we conclude that
$$
\frac{1}{\log X} \Big|\sum_{n\le X} \frac{f(n)}{n} \Big|
\le \frac 12 +\frac{1}{\log u} \int_{\sqrt{u}}^{u} (1-I(t))\frac{dt}{t}
+o(1) \le \frac 43 - \frac{1}{e^{2/3}} +o(1) .
$$

\subhead 3c.  Proof of Theorem 1 for $k=4$ \endsubhead

\noindent Note that $g(p)=\sqrt 2 -1$ if $f(p)=\pm i$, and $g(p)=f(p)$
if $f(p)=0$, $\pm 1$.  We may assume that $I(X)\le 1$ else, applying Lemmas
3 and 4 (taking $u$ to be the largest integer with $I(u)\le 1$),
we deduce that $|\sum_{n\le X} f(n)/n| \le (2-2/\sqrt{e}+o(1))\log X$.

Let
$$
A = \sum\Sb p\le X \\ f(p)=0, \pm i\endSb \frac{1}{p}, \qquad
\text{and }\qquad B= \sum\Sb p\le X\\ f(p)=-1\endSb \frac 1p.
$$
so that $1\ge I(X) \ge (2-\sqrt{2})A+2B$. For $t\le X$ we have
$$
\align
I(t) & \ge (2-\sqrt{2}) \sum\Sb p\le t\\ f(p)=0, \pm i\endSb \frac 1p +
2\sum\Sb p\le t \\ f(p)=-1 \endSb \frac 1p \\
& \geq (2-\sqrt{2}) \biggl(
A - \sum\Sb t<p\le X\\ f(p)=0, \pm i\endSb \frac 1p\biggr) +
2\biggl( B-\sum\Sb p\le t \\ f(p)=-1 \endSb \frac 1p \biggr) \\
&\ge
\cases
(2-\sqrt{2})A+ 2(B-\log(\log X/\log t)) +o(1) &\text{if   } X^{e^{-B}} \le
t\le X\\
(2-\sqrt 2)(A+ B -\log (\log X/\log t))+o(1) &\text{if   } X^{e^{-A-B}}
\le t\le X^{e^{-B}}.\\
\endcases \\
\endalign
$$
Of course $I(t)\geq 0$ for $t\le X^{e^{-A-B}}$.
Using these lower bounds and the third bound in Lemma 4 we deduce that
$$
\align
\frac{1}{\log X} \left| \sum_{n\le X} \frac{f(n)}{n} \right|
&\le
 1 - \frac{1}{\log X}\int_{\sqrt{X}}^{X} I(t) \frac{dt}{t}+o(1)\\
&\le 1 -
\frac{\sqrt 2}{\log X}
\int_{\max\{ X^{e^{-B}},X^{1/2}\} }^X \left(B-\log\left(\frac{\log X}{\log t}
\right) \right)\frac {dt}{t}\\
&\hskip .3 in - \frac{(2-\sqrt{2})}{\log X}\int_{\max\{ X^{e^{-A-B}},X^{1/2}\} }^{X}
\left(A+B-\log \left(\frac{\log X}{\log t}\right)\right) \frac{dt}{t} +o(1)
\\
&\leq F(A,B)+o(1) ,\\
\endalign
$$
say, where (using (2.2) to compute the integrals)
$$
F(A,B)=\cases
3-2B-(2-\sqrt{2})(A+e^{-A-B})-\sqrt{2}e^{-B}
& \text{\rm if}\ A+B\leq \log 2 \\
2+1/\sqrt{2} - (1-1/\sqrt{2})(A+\log 2)-\sqrt{2}e^{-B}- (1+1/\sqrt{2})B
 & \text{\rm if}\ B\leq \log 2 \leq A+B\\
2 - \log 2 -B - (1-1/\sqrt{2})A & \text{\rm if}\ B\geq \log 2.\\
\endcases
$$
By differentiation we find that $F(A,B)$ is a non-increasing function of both
$A$ and $B$, for $A,B\geq 0$.

Since Re$ f(n)\ge 1- \sum_{p|n} (1-\text{Re } f(p))$ we have, by (1.2),
$$
o(X) = \sum_{n\le X}\text{Re } f(n) \ge X -\sum_{p\le X} (1-\text{Re }f(p))
\frac{X}{p} + o(X),
$$
so that $A+2B = \sum_{p\le X} (1-\text{Re }f(p))/p \ge 1+o(1)$.
Therefore
$$
F(A,B) \le \min_{0\leq A\leq 1} F(A,(1-A)/2) +o(1)
= F(A_0,(1-A_0)/2) +o(1) = 0.8296539741\ldots,
$$
where $A_0:=2\log( (3-\sqrt{2})/2)+1\approx 0.5358665582\ldots$.

\head 4.  Proof of Theorem 1 for $S=\Bbb U$\endhead

\noindent In this proof we will assume only that (1.2) holds for all $x$ in
the interval $X\leq x\leq X^{\sqrt{e}}$.

We may assume that  $I(X)<1$ else by applying Lemmas 3 and 4
(as before) we get that
$(1/\log X)  |\sum_{n\le X} f(n)/n| \le (2-2/\sqrt{e}) + o(1)$.
We may also assume that $I(\sqrt X)\le \frac 2{35}$ else, by Lemma 4,
$\frac{1}{\log X} |\sum_{n\le X} f(n)/n| \le 1-I(\sqrt{X})/2 +o(1) \le
\frac{34}{35}+o(1)$.

Define
$$
A=\sum_{p\le X} \frac{1-\text{Re }f(p)}p
\qquad \text{and} \qquad
B=\sum_{p\le X} \frac{1-\text{Re }f(p)}p \frac{\log (X/p)}{\log X}.
$$
The second bound in Lemma 4 gives that
$$
\align
\frac{1}{\log X} \Big|\sum_{n\le X} \frac{f(n)}{n} \Big|
&\le 1-\frac{1-I(\sqrt{X})}{\log X} \int_{1}^{X} I(t) \frac {dt}{t} +o(1)
\le 1 - \frac{33}{35\log X} \int_1^X I(t) dt + o(1)\\
&= 1- \frac{33}{35} \sum_{p\le X} \frac{1-g(p)}{p} \frac{\log (X/p)}{\log X}
+o(1) \leq 1 -\frac{33B}{70} + o(1) , \tag{4.1}  \\
\endalign
$$
where the final inequality holds since $1-g(p) = 2-|1+f(p)|
\ge (1-\text{Re } f(p))/2$.

Now for $t\le X$ we
have
$$
\sum_{p\le t} \frac{1-\text{Re } f(p)}{p}
= A-\sum_{t\le p\le X} \frac{1-\text{Re } f(p)}{p}
\ge A- 2\log \left(\frac{\log X}{\log t}\right) +o(1),
$$
and also note that $\sum_{p\le t} \frac{1-\text{Re } f(p)}{p}  \ge 0$.
Hence, using (2.2), we deduce that
$$
\align
B&=\frac{1}{\log X} \int_{1}^{X} \sum_{p\le t} \frac{1-\text{Re }f(p)}{p}
\frac{dt}{t}
\ge \frac{1}{\log X} \int_{X^{e^{-A/2}}}^X
\left( A- 2\log \left(\frac{\log X}{\log t}\right) +o(1)\right) \frac{dt}{t}\\
&= A-2+2e^{-A/2} +o(1)\geq A_0-2+2e^{-A_0/2} +o(1), \\
\endalign
$$
for all $A\geq A_0:=2\log (2(\sqrt{e}-1))=0.5207901030\ldots$ (since $A-2+2e^{-A/2}$ increases
for all $A\ge 0$).

Put $h(n) =\sum_{ab=n} f(a) \overline{f(b)}$, so that $h(n) \ge d(n)
(1-\sum_{p^k|n} (1-\text{Re} f(p)))$. Therefore
$$
\sum_{n\le Y} h(n) \ge Y\log Y - 2Y \sum_{p\le Y}
\frac{1-\text{Re }f(p)}{p} \log \frac{Y}{p} +o(Y\log Y).
$$
Combining this with Lemma 5 we deduce that for $X\le Y\le X^{\sqrt{e}}$
$$
\sum_{p\le Y} \frac{1-\text{Re }f(p)}{p} \log \frac Yp
\ge \log \frac{Y}{X} +o(\log X).
$$
Taking $Y=X^{1+\alpha}$ with $0\leq \alpha \leq \sqrt{e}$, and using that $1-\text{Re }f(p)\leq 2$,
we deduce that
$$
\align
\alpha+o(1)&\leq \alpha A + B + 2\sum_{X\le p\le X^{1+\alpha}}
\frac{\log (X^{1+\alpha}/p)}{p\log X} \\
&= \alpha A + B + 2(1+\alpha)\log (1+\alpha)-2\alpha +o(1).\\
\endalign
$$
If $A\le 1$ then taking $\alpha = e^{(1-A)/2} -1$ we deduce from the
above that
$B\ge 2 e^{(1-A)/2} +A -3 +o(1) \ge 2e^{(1-A_0)/2} +A_0-3$,
for $0\leq A \leq A_0
$ (since the function here is decreasing for $0\leq A\leq 1$).

Either way we deduce that $B\ge 0.062284\ldots$ so that
$1-33B/70\le \frac{34}{35}$
which, by (4.1),  proves the desired estimate for $\gamma({\Bbb U})$.

\remark{Remark} The constant $34/35=.9714285714\ldots$ may be replaced
by $.9706838406\ldots$ in the above proof.
\endremark

\head 5.  Upper bounds on average:  Proof of Theorem 3 \endhead

\noindent Let $k\ge 2$ be an integer, and suppose that
$f\in {\Cal F}(S_k)$ satisfies (1.3) for $X\le x\le X^{\vphi(k)+1}$.
As in section 2 we can assume that $f(p)=1$ for all $p\leq y$, without loss
of generality. Let $g(n)$ be the multiplicative function defined by
$$
g(n) = \sum_{ \prod_{(j,k)=1} a_j = n} \prod\Sb j=1\\(j,k)=1\endSb^k
f(a_j)^j,
$$
and $h(n) =\sum_{d|n} g(d)$. We consider
$$
\sum_{n\le X^{\vphi(k)+1}} h(n)
= \sum_{a_0\prod_{(j,k)=1} a_j\le X^{\vphi(k)+1}} \prod_{(j,k)=1}
f(a_j)^j,
$$
and distinguish two types of terms: when all the $a_j$ with $(j,k)=1$
are below $X$, and when one of them exceeds $X$.  The first type
contribute
$$
\align
&\sum\Sb a_j\le X \\ (j,k)=1\endSb \prod_{(j,k)=1} f(a_j)^j
\biggl( X^{\vphi(k)+1}\prod_{(j,k)=1} \frac{1}{a_j } + O(1)\biggr)
\\
= &X^{\vphi(k)+1} \prod_{(j,k)=1} \biggl(  \sum_{n\le X} \frac{f(n)^j}{n}
\biggr)
+ O(X^{\vphi(k)})
\\
\endalign
$$
We next consider the contribution of the second type of terms.  Suppose
for example that $a_1 > X$ is the largest of the $a_j$'s with $(j,k)=1$.
For fixed $a_0$, $a_j$ ($j\ge 2$ with
$(j,k)=1$), we get,
from our assumption on the range in which (1.3) holds,
that the sum over $a_1$ is
$$
o\biggl( \frac{X^{\vphi(k)+1}}{a_0} \prod\Sb j\ge 2\\ (j,k)=1\endSb
\frac{1}{a_j}\biggr).
$$
Writing $n=a_0 \prod_{j\ge 2, (j,k)=1} a_j$ we see that
the contribution of these terms is
$$
\align
o\biggl( X^{\vphi(k)+1}\sum_{n\le X^{\vphi(k)}}
\frac{d_{\vphi(k)}(n)}{n}\biggr)
&= o\left(X^{\vphi(k)+1} e^{\vphi(k)} \sum_{n=1}^{\infty}
\frac{d_{\vphi(k)}(n)}{n^{1+1/\log X}} \right)\\
&= o\left(X^{\vphi(k)+1} e^{\vphi(k)}
\zeta\left(1+\frac 1{\log X}\right)^{\vphi(k)}\right)
\\
&= o( X^{\vphi(k)+1} (\log X)^{\vphi(k)}).
\tag{5.1}\\
\endalign
$$
The same argument applies when any other $a_j$ is the largest.  Thus we
conclude that
$$
\frac{1}{X^{\vphi(k)+1}} \sum_{n\le X^{\vphi(k)+1}} h(n)
=  \prod_{(j,k)=1} \sum_{n\le X} \frac{f(n)^j}{n} + o\left((\log X)^{\vphi(k)}\right).
\tag{5.2}
$$

Since $h(n)=\sum_{d|n} g(d)$ and $|g(n)| \le d_{\vphi(k)}(n)$
we have
$$
\frac{1}{X^{\vphi(k)+1}} \sum_{n\le X^{\vphi(k)+1}} h(n)
= \sum_{d\le X^{\vphi(k)+1}} \frac{g(d)}{d} + O\biggl(
\frac{1}{X^{\vphi(k)+1}} \sum_{n\le X^{\vphi(k)+1}} d_{\vphi(k)}(n)\biggr).
$$
Writing $d_{\vphi(k)}(n) =\sum_{a|n} d_{\vphi(k)-1}(a)$, and
arguing as in (5.1),  we see that
$$
\frac{1}{X^{\vphi(k)+1}} \sum_{n\le X^{\vphi(k)+1}} d_{\vphi(k)}(n)
\le \sum_{a\le X^{\vphi(k)+1}} \frac{d_{\vphi(k)-1}(a)}{a}
\ll (\log X)^{\vphi(k)-1}.
$$
These observations and (5.2) give that
$$
\align
\prod_{(j,k)=1}  \biggl| \sum_{n\le X} \frac{f(n)^j}{n}\biggr|  &=
\Big|\sum_{d\le X^{\vphi(k)+1}} \frac{g(d)}{d}\Big|
+ o\left((\log X)^{\vphi(k)}\right) \\
&\le \sum_{d\le X^{\vphi(k)+1}} \frac{|g(d)|}{d}  +
o\left((\log X)^{\vphi(k)}\right).\\
\endalign
$$
Write $\delta = c/\log X$ for some positive constant $c>0$ to be fixed later.
Then
$$
\align
\sum_{d\le X^{\vphi(k)+1}} \frac{|g(d)|}{d} &\le e^{c(\vphi(k)+1)}
\sum_{d=1}^{\infty} \frac{|g(d)|}{d^{1+\delta}} \\
&\ll e^{c\vphi(k)} \zeta(1+\delta)^{\vphi(k)} \prod_p \left\{ 1 + \frac{|g(p)|}{p^{1+\delta}}
+ \frac{|g(p^2)|}{(p^2)^{1+\delta}} + \dots \right\} \left( 1 - \frac 1{p^{1+\delta}} \right)^{\vphi(k)} \\
&\ll \left\{  \frac{e^c}c  \log X  \exp\left(- \sum_p\frac{ 1-|g(p)|/\vphi(k)}{p^{1+\delta}} \right)\right\}^{\vphi(k)}.
\\
\endalign
$$
To justify this last step note
that the $p$th term in the Euler product is $1$ when $f(p)=1$, which
happens for all primes $p\leq y$,
and so the error term for the whole Euler product, in the transition
from the penultimate bound to the last one, is
$\prod_{p>y} \exp ( O( \phi(k)^2/p^2)) = 1+o(1)$.
We deduce that
$$
\left(\prod_{(j,k)=1} \frac{1}{\log X} \left|
\sum_{n\le X} \frac{f(n)^j}{n}\right| \right)^{\frac{1}{\vphi(k)}}
\leq  \frac{e^c}c
\exp\left(- \sum_p\frac{ 1-|g(p)|/\vphi(k)}{p^{1+\delta}} \right)+
o_k(1).
\tag{5.3}
$$

For each prime $p\nmid q$ let $l_p$ be such that $f(p)$ is a
primitive $l_p$-th root of unity.  Note that $l_p$ is a
divisor of $k$, and that as $j$ varies over all reduced
residue classes $\pmod k$, $f(p)^j$ runs over all
primitive $l_p$-th roots of unity $\vphi(k)/\vphi(l_p)$ times.
Thus
$$
g(p) = \sum_{(j,k)=1} f(p)^j
= (\vphi(k)/\vphi(l_p)) \sum\Sb a\pmod{l_p}\\(a,l_p)=1\endSb e(a/l_p)
= \mu(l_p) \vphi(k)/\vphi(l_p).
$$
 Define $l_p=k$ if $p|q$, and note
that here $g(p)=0$.  From these remarks and (5.3) it follows
that
$$
\left(\prod_{(j,k)=1} \frac{1}{\log X} \left|
\sum_{n\le X} \frac{f(n)^j}{n}\right| \right)^{\frac{1}{\vphi(k)}}
\leq \frac{e^c}c
\exp\left(-  \sum_{p\le X} \frac{1-1/\vphi(l_p)}{p^{1+\delta}}\right)+o_k(1).
\tag{5.4}
$$

To estimate the right hand side of (5.4) we employ the
following result of Hildebrand [5] together
with an idea of Vinogradov as exploited by Norton [6].

\proclaim{Lemma 6}  Fix $\theta>0$. We have
$$
\lim_{x \to \infty} \inf \frac 1x \sum\Sb n\le x \\ (n,P)=1 \endSb 1
=  \rho(e^{\theta}) ,
$$
where $\inf$ is taken over all subsets $P$ of the primes up to $x$,
such that  $\sum_{p\in P} 1/p = \theta + o(1)$.
Here $\rho(u)$ is the Dickman-de Bruijn function, defined by $\rho(u) =1$
for $0 \le u \le 1$, and $u\rho^{\prime}(u) = -\rho(u-1)$ for all
$u \ge 1$.
The lower bound is attained when $P$ is the set of primes in
$[x^{e^{-\theta}},x]$.
\endproclaim

Let $l$ be a divisor of $k$ with $l<k$, and let $P_{l}$ denote
the product of those
primes $p$ below $X$ for which $l_p \nmid l$.  Observe that if $p\nmid P_l$
then $f^l(p)=1$.   Note that
$$
\sum\Sb n\le X \\ f^l(n)=1\endSb 1 = \sum_{n\le X} \frac{l}{k}
\sum_{v=1}^{k/l} f^{lv}(n)
= \left(\frac{l}{k}+o(1)\right) X,
$$
using (1.4) when $1\le v <k/l$.  On the other hand
$$
\sum\Sb n\le X \\ f^l(n)=1\endSb 1 \ge \sum\Sb n\le X \\
(n,P_l)=1\endSb  1 \geq  \rho\Big(   \exp\Big(
\sum_{p|P_l}  \frac 1p   \Big)\Big)  X +o(X)
$$
by Lemma 6, so that
$$
\exp \Big( \sum\Sb p\le X \\ l_p \nmid l\endSb \frac 1p \Big)
\ge \rho^{-1}(l/k) :=\theta^{-1} \sim
\frac{\log k}{\log\log k} \quad \text{\rm if}\ l = k^{o(1)},
\tag{5.5}
$$
since $\rho(u)=e^{-u\log u (1+o(1))}$.
Here given $x\in [0,1)$, $\rho^{-1}(x)$ denotes the unique $u$ with
$\rho(u)=x$.

If $k=\prod p^{\alpha_p}$ then take $l= \prod p^{\beta_p}$
where $\beta_p=[\min \{ \alpha_p,\log \log k/2\log p\} ]$.  Note that
$p^{\beta_p}\leq \sqrt{\log k}$ so that
$l\le \prod_{p\le \sqrt{\log k}} \sqrt{\log k} \le \exp( (\log k)^{\frac 23})$.
Further if $l_p \nmid l$ then $l_p$ is divisible by a
prime power larger than $\sqrt{\log k}$ (as $l_p$ divides $k$), and so
$\phi(l_p) \ge \sqrt {\log k}/2$.
   From these remarks and (5.5) we see that
$$
\align
\exp\Big(-\sum_{p\le X} \frac{1-1/\vphi(l_p)}{p^{1+\delta}} \Big)
&\le \exp\Big(-
(1+O((\log k)^{-\frac 12}))
\sum\Sb p\leq X\\ l_p\nmid l \endSb \frac 1{p^{1+\delta}} \Big) \\
&\le \{ 1+o_k(1)\} \exp\Big(- \sum\Sb X^\theta<p\leq X\endSb
\frac 1{p^{1+\delta}} \Big) \\
&\le \{ 1+o_k(1)\}
\frac{\log\log k}{\log k}
\exp \left( \int^1_0 \left( \frac{1-e^{-ct}}t   \right) dt\right). \\
\endalign
$$
Combining this with (5.4) gives
$$
\Big(\prod_{(j,k)=1} \frac{1}{\log X} \Big|
\sum_{n\le X} \frac{f(n)^j}{n}\Big| \Big)^{\frac{1}{\vphi(k)}}
\le \{ 1+o_k(1)\}
\frac{\log\log k}{\log k}
\exp \Big( c-\log c +\int^1_0 \Big( \frac{1-e^{-ct}}t   \Big) dt\Big).
$$
The left side is minimized at $c=0.5671432904\ldots$
(that is, where $e^{-c}=c$), giving an upper bound
$<  2.8661 e^\gamma \log \log k/\log k<(43/15)e^\gamma \log\log k/\log k$,
proving Theorem 3.

\head 6.  Integral equations:  Proofs of Theorems 2a and 2b \endhead

\noindent Proposition 1 of [3] shows how problems concerning the
distribution of multiplicative
functions (of absolute values $\leq 1$), and problems
concerning certain integral equations
are essentially equivalent.
In particular, our questions on $\gamma(S)$ and
$\gamma_k$ may be reformulated in terms of integral equations.
Regarding the proofs of Theorems 1 and 3, this confers only a marginal
advantage and so we did not pursue this approach in those contexts.
However the integral equations approach considerably simplifies the
treatment of the lower bounds for $\gamma(S_k)$ and $\gamma_k$
claimed in  Theorems 2a and 2b.  We begin by recapitulating the
relevant material from [3].

For a given closed, subset $S$ of the unit disc, let $K(S)$ denote
the class of measurable functions $\chi: [0,\infty) \to S^*$ (the convex hull
of $S$)
with $\chi(t)= 1$ for $0\le t\le 1$. There is a unique (continuous)
$\sigma: [0,\infty) \to {\Bbb U}$ satisfying
$$
\align
u\sigma(u)=  \int_0^u \sigma(u-t)\chi(t) dt
 \ \ &\text{\rm for} \ \ u > 1 ,\tag{6.1} \\
\text{\rm with the initial condition} \  \ \sigma(u)=1 \ \ \
&\text{\rm for} \ \ 0\le u \le 1. \\
\endalign
$$
Define $\Lambda(S)$ to be the set of such values $\sigma(u)$.

\proclaim{Proposition 1} Let $f$ be a multiplicative function with $|f(n)|
\le 1$ for all $n$, and $f(n)=1$ for
$n\le y$.  Let $\vartheta(x) =\sum_{p\le x} \log p$ and define
$$
\chi(u) = \chi_f(u) = \frac{1}{\vartheta(y^u) } \sum_{p\le y^u} f(p)\log p.
$$
Then $\chi(t)$ is a measurable function taking values in the unit disc
with $\chi(t)=1$ for $t\le 1$, and
$\sigma(u)$, the corresponding unique solution to {\rm (6.1)}, satisfies
$$
\frac{1}{y^u} \sum_{n\le y^u} f(n) =\sigma(u) +O\left(\frac{u}{\log y}\right).
$$
\endproclaim

The converse to Proposition 1 is also true.

\proclaim{Proposition 1 (Converse)} Let $S\subset {\Bbb {\Bbb U}}$ and
$\chi \in K(S)$ be given.  Given $\epsilon >0$ and $u\ge 1$ there
exist arbitrarily large $y$ and $f\in {\Cal F}(S)$ with $f(n)=1$ for
$n\le y$ and
$$
\left| \chi(t) - \frac{1}{\vartheta(y^t)} \sum_{p\le y^t} f(p)\log p
\right| \le \epsilon \ \ \text{for almost all } 0 \le t\le u.
$$
If $\sigma(u)$ is the solution to {\rm (6.1)} for this $\chi$ then
$$
\sigma(t) = \frac{1}{y^t} \sum_{n\le y^t} f(n) + O(u^{\epsilon}- 1)
+O\left(\frac{u}{\log y}\right) \ \ \text{for all } t\le u.
$$
\endproclaim




Theorems 2a and 2b will follow from the following result on integral equations.

\proclaim{Proposition 2} For $0\leq \delta \leq 1$
define $\chi_\delta(t)=1$ for $0\leq t\leq 1$,
and $\chi_\delta(t)=-\delta$ for $t\geq 1$.
Let $\sigma_\delta$ denote the corresponding solution in {\rm (6.1)}.
We have $\sigma_\delta(u)= 1-(1+\delta)\log u$ for $1\le u\le 2$.
\item{(i)} For $0<\delta\leq 1$ there exists a positive real
root of $\sigma_\delta(u)=0$. If $U_\delta$
is the smallest such root then $U_\delta$ is a decreasing function of
$\delta$ with $U_\delta=e^{1/(1+\delta)}$ when $1/\log 2 -1 \le \delta
\le 1$,
and $U_\delta \sim
\log(1/e^{\gamma}\delta)/\log\log(1/e^{\gamma}\delta)$ as $\delta \to 0$.
\item{(ii)} There exists $\chi$ with $\chi(t)=-\delta$
for $1\leq t\leq U_\delta$ and $\chi(t)\in [-\delta,1]$ for $t>U_\delta$,
such that $\sigma(u)=\sigma_\delta(u)$ for $0\le u \le U_\delta$, and
$\sigma(u)=0$
when $u\geq U_\delta$.  The function $\chi(t)$ is
continuous for all $t>U_\delta$.  As $\delta \to 0^+$ we have that
$$
I_\delta:= \frac 1{U_\delta} \int_0^{U_\delta} \sigma(t) dt =
\frac{e^\gamma}{U_\delta} +O(\delta).
$$
\endproclaim

Assuming Proposition 2 for the moment, we complete the
proofs of Theorems 2a, b.

\demo{Proof of Theorem 2a}  Taking $\delta=1$ in Proposition 2 we know
that there is a $\chi \in K(S_2)$ such that $\sigma(u) =1$ for $u\le 1$,
$\sigma(u)=1-2\log u$ for $1\le u\le \sqrt{e}$ and
$\sigma(u)=0$ for $u\ge \sqrt{e}$.  By Proposition 1 (Converse) we
can find $f \in {\Cal F}(S_2)$ such that for $0\le t\le A\sqrt{e}$,
$$
\frac{1}{y^t} \sum_{n\le y^t} f(n) = \sigma(t) + O(\epsilon A + A/\log y).
$$
Thus with $X=y^{\sqrt{e}}$ we find that $f$ satisfies (1.2) in
the range $X\le x\le X^A$, and further
$$
\frac{1}{\log X}  \sum_{n\le X} \frac{f(n)}{n}
= \frac{1}{\sqrt{e}} \int_0^{\sqrt{e}} \frac{1}{y^t} \sum_{n\le y^t} f(n)
dt + o(1) = 2-\frac{2}{\sqrt{e}} + O(\epsilon A +A/\log y)+o(1).
$$
   From this it follows that $\gamma(S_2;A) \ge 2-2/\sqrt{e}$, and
since $A$ is arbitrary the same holds for $\gamma(S_2)$.

If $k\ge 3$ is odd then we apply Proposition 2 with
$\delta= \delta_k (\ge 1/2)$.  We thus find
$\chi \in K([-\delta_k,1]) \subset
K(S_k)$ such that $\sigma(u) =1-(1+\delta_k) \log u$ for
$1\le u\le U_{\delta}= e^{1/(1+\delta_k)}$, and $\sigma(u)=0$ for
$u\ge U_\delta$.  We now argue exactly as above, constructing
$f$ via Proposition 1 (converse), and noting that
$(1/U_\delta) \int_0^{U_\delta} \sigma(t) dt =
(1+\delta)(1- e^{-1/(1+\delta)})$.  This proves Theorem 2a.

\enddemo

\demo{Proof of Theorem 2b} Here we apply Proposition 2 with $\delta
= 1/(k-1)$.  So there is a $\chi \in K([-1/(k-1),1])$ such that
$\sigma(u) = 0$ for $u\ge U_\delta$.  Write $\chi(t) = 1- k\alpha(t)$
so that $0\le \alpha(t)\le 1/(k-1)$ for all $t$; and note that
$\alpha(t)$ is piecewise linear for $t\leq U_\delta$ and continuous
for $t>U_\delta$.  For fixed $A\ge 1$, $\epsilon >0$
and large $y$ we may easily partition
the set of primes below $y^{AU_\delta}$  as
${\Cal P}_0\cup {\Cal P}_1 \cup \ldots \cup {\Cal P}_{k-1}$
such that for $0\le t\le AU_\delta$ we have
$$
\frac{1}{y^t} \sum\Sb p\le y^t\\ p\in {\Cal P}_0 \endSb \log p
= 1- (k-1)\alpha(t) +O(\epsilon),
$$
and for $j=1$, $2$, $\ldots$, $k-1$ we have
$$
\frac{1}{y^t} \sum\Sb p\le y^t \\ p\in {\Cal P}_j\endSb \log p
= \alpha(t) + O(\epsilon).
$$
We now choose $f \in {\Cal F}(S_k)$ by setting
$f(p)=e(\ell/k)$ if $p\in {\Cal P}_\ell$ for $\ell =0,\ldots, k$.
For such $f$ we get that for $j=1,\ldots,k-1$, and $t\le AU_\delta$,
$$
\frac{1}{y^t} \sum_{p\le y^t} f(p)^j \log p = 1- k\alpha(t) + O(k\epsilon)
= \chi(t) +O(k\epsilon),
$$
so that by Proposition 1 (Converse) we may conclude that
$$
\frac{1}{y^t} \sum_{n\le y^t} f(n)^j = \sigma(t) + O(AkU_\delta \epsilon
+ AU_\delta/\log y).
$$
Thus with $X=y^{U_\delta}$ we find that (1.3) holds in the range
$X\le x\le X^A$, and further
that
$$
\frac{1}{\log X} \sum_{n\le X} \frac{f(n)^j}{n}
= \frac{1}{U_\delta}\int_0^{U_\delta} \sigma(t) dt + o(1) = I_\delta +o(1).
$$
This gives $\gamma_k(A) \ge I_{\delta}$ and, since $A$ is
arbitrary, Theorem 2b follows.

\enddemo

To prove Proposition 2 we require Lemma 3.4 from [3] which we
quote below.

\proclaim{Lemma 7}  Let $\chi$ and ${\hat \chi}$ be two elements of
$K({\Bbb U})$, and let $\sigma$ and $\hsigma$
be the corresponding solutions to {\rm (6.1)}.  Then $\sigma(u)$ equals
$$
\hsigma(u) + \sum_{j=1}^{\infty} \frac{(-1)^j}{j!} \int\Sb
t_1,\ldots ,t_j \ge 1\\ t_1+\ldots +t_j \le u \endSb
\frac{\hchi(t_1)-\chi(t_1)}{t_1}
\ldots \frac{\hchi(t_j)-\chi(t_j)}{t_j}
\hsigma(u-t_1-\ldots-t_j)dt_1\ldots dt_j.
$$
\endproclaim

\demo{Proof of Proposition 2}  First we apply Lemma 7 taking
$\hchi(t)=1$ for $t\le 1$
and $\hchi(t)=0$ for $t>1$ (so that $\hsigma(u)=\rho(u)$), and $\chi(t)=1$
for all $t\ge 0$ (so that $\sigma(u)=1$).  This gives that
$$
1=\rho(u) + \sum_{j=1}^{\infty} \frac{1}{j!} \int\Sb t_1,\ldots, t_j\ge 1\\
t+1+\ldots+t_j \le u\endSb \rho(u-t_1-\ldots -t_j)
\frac{dt_1 dt_2 \ldots dt_j}{t_1t_2\ldots t_j},
$$
and thus
$$
0 \le \sum_{j=2}^{\infty} \frac{1}{j!} \int\Sb t_1,\ldots, t_j\ge 1\\
t+1+\ldots+t_j \le u\endSb \rho(u-t_1-\ldots -t_j)
\frac{dt_1 dt_2 \ldots dt_j}{t_1t_2\ldots t_j} \le 1. \tag{6.2}
$$

Next we apply Lemma 7 taking $\hchi$ as above, and $\chi(t) =\chi_\delta(t)$.
Hence
$$
\sigma_{\delta}(u)
=\rho(u) + \sum_{j=1}^{\infty} \frac{(-\delta)^j}{j!} \int\Sb
t_1,\ldots ,t_j \ge 1\\ t_1+\ldots +t_j \le u \endSb
\ \rho(u-t_1-\ldots-t_j) \frac{dt_1 dt_2\ldots dt_j}{t_1 t_2\ldots t_j},
\tag{6.3}
$$
and in view of (6.2) we may conclude that
$$
\left| \sigma_{\delta}(u)
-\left( \rho(u) -\delta \int_1^u \rho(u-t) \frac{dt}t
\right) \right|\leq \delta^2. \tag{6.4}
$$
If $\delta$ is sufficiently small then the asymptotics
$\rho(u)=e^{-u\log u (1+o(1))}$, and $\int_1^u \rho(u-t) dt/t
\sim (1/u) \int_0^{\infty} \rho(t) dt \sim e^{\gamma}/u$
enable us to deduce that $U_\delta \sim \log(1/e^{\gamma} \delta)/\log \log
(1/e^{\gamma} \delta)$.  Thus for sufficiently small $\delta$
we have established the existence of a positive real root of
$\sigma_\delta(u)=0$, and the asymptotic for $U_\delta$ claimed in
part (i).

Suppose as above that $\delta$ is sufficiently small so that
$U_\delta$ exists, and
let $1\geq \eta \ge \delta  >0$.  Take $\chi=\chi_\delta$ and $\hchi =
\chi_\eta$ in Lemma 7.  Evaluating at $u=U_\delta$ we conclude that
$$
0=\sigma_\delta(U_\delta) = \sigma_\eta(U_\delta) +
\sum_{j=1}^{\infty} \frac{(\eta-\delta)^j}{j!}
\int\Sb t_1,\ldots, t_j \ge 1\\ t_1 + \ldots+ t_j\le u\endSb
\sigma_\eta(u-t_1-\ldots -t_j) \frac{dt_1 dt_2 \ldots dt_j}{
t_1 t_2 \ldots t_j}.
$$
It follows at once that $\sigma_\eta(u)$ must change sign (and hence have
a zero) in $(0,U_\delta]$.  Thus $U_\eta$ exists for all $\eta \in (0,1]$
and is a decreasing function of $\eta$.

Lastly note that for $1\le u \le 2$ we
have $\sigma_{\delta}(u) = 1- (1+\delta)\log u$, from which it follows
that in the range $(0.4426\ldots =)1/\log 2 - 1
\le \delta \le 1$ we have $U_\delta = e^{1/(1+\delta)}$.
This completes the proof of part (i).

We now turn to the proof of part (ii).  First observe that
$\sigma_\delta(u)$ satisfies the differential-difference equation
$u\sigma_{\delta}^{\prime}(u) = - (1+\delta) \sigma_{\delta}(u-1)$, for
$u\ge 1$.  Since $\sigma_{\delta}(u)$ is positive for $0\le u< U_\delta$
we conclude that $\sigma_{\delta}^{\prime}(u) <0$ for $1\le u
< U_\delta+1$.  Further observe that
$$
1= \sigma_\delta(1)-\sigma_\delta(U_\delta) = \int_1^{U_\delta}
(-\sigma_\delta^{\prime}(t) ) dt.
\tag{6.5}
$$

We now take $\chi(t) = 1$ for $t\le 1$, $\chi(t) = -\delta$ for
$1\le t\le U_\delta$, and for $t> U_\delta$ we define
$$
\chi(t) = \int_1^{U_\delta} (-\sigma_\delta^{\prime}(v)) \chi(t-v) dv.
$$
   From this definition it is clear that $\chi$ is continuous for $t>U_\delta$.
We shall first show that $\chi(t) \in [-\delta,1]$ for all $t$.  Plainly
this holds for all $t\le U_\delta$.  From our definition of $\chi$,
the positivity of $-\sigma_{\delta}^{\prime}(t)$ in $(1,U_{\delta})$,
and (6.5) we immediately glean that
$$
\chi(t) \le \max_{z\in [t-1,t-U_\delta]} \chi(z),
\qquad \text{and} \qquad \chi(t) \ge \min_{z \in [t-1,t-U_\delta]}
\chi(z).
$$
Inductively it follows that $\chi(t) \in [-\delta,1]$ for all $t$
as desired.

Next we demonstrate that if $\sigma(u)$ denotes the solution
to (6.1) for $\chi$ constructed above, then $\sigma(u)=0$ for
$u\ge U_\delta$.  Note that for $u\le U_\delta$ we have
$\sigma(u)=\sigma_\delta(u)$, and so in particular
$\sigma(U_\delta) =\sigma_\delta(U_\delta)=0$.
Now, since $\sigma'(t)=0$ for $0<t<1$, we have for $u\geq U_\delta$,
$$
\frac d{du} \int^u_{u-U_\delta} \sigma(u-t)\chi(t) dt
=
\int^u_{u-U_\delta} \sigma'(u-t)\chi(t) dt
+ \chi(u) \sigma(0) - \sigma(U_\delta) \chi(u-U_\delta) =0
$$
by definition of $\chi(u)$. Therefore
$\int^u_{u-U_\delta} \sigma(u-t)\chi(t) dt$ is a constant for
$u\geq U_\delta$, and at $u=U_\delta$ equals, by definition,
$U_\delta \sigma (U_\delta)=0$.  Hence for $u\geq U_\delta$,
$$
u\sigma(u) = \int^u_0 \sigma(u-t)\chi(t) dt
= \int^{u-U_\delta}_0 \sigma(u-t)\chi(t) dt
= \int^u_{U_\delta} \sigma(v)\chi(u-v) dv .
$$
We claim that this gives $\sigma(u)=0$ for all $u\geq U_\delta$.
If not, select $u>U_\delta$ such that $|\sigma(u)| >0$ and
such that $|\sigma(u)| \geq |\sigma(v)|$ for all
$v\in [U_\delta,u]$; then
$$
u|\sigma(u)| \leq  \int^u_{U_\delta} |\sigma(v)\chi(u-v)| dv
\leq |\sigma(u)| \int^u_{U_\delta}  dv =(u-U_\delta) |\sigma(u)|,
$$
giving a contradiction.

We have thus constructed $\chi$ and $\sigma$
as desired.  For $u\le U_\delta$ we have $\sigma(u)=\sigma_\delta(u)
= \rho(u)+O(\delta /u +\delta^2)$ by (6.4),
and so
$$
\frac 1{U_\delta} \int_0^{U_\delta} \sigma(t) dt =
\frac 1{U_\delta} \int_0^{U_\delta} \Big(\rho(t) +O\Big(\frac{\delta }{(t+1)}
+\delta^2\Big)\Big) dt
=\frac {e^\gamma}{U_\delta} +O(\rho(U_\delta) +\delta),
$$
which completes the proof of Proposition 2.

\enddemo

We conclude this section with some numerical data
pertaining to Proposition 2.
We noted earlier that in the range $0.44269\ldots = 1/\log 2 - 1\le
\delta \le 1$ we have $U_\delta =e^{1/(1+\delta)}$
lying in the range $[1,2]$.When $2\leq u\leq 3$ we have
$$
\sigma_\delta(u)=1-(1+\delta)\log u
+ \frac{(1+\delta)^2}2 \int_1^{u-1} \frac{\log(u-t)}t dt.
$$
   From this we find that $U_{.061129446\ldots}=3$,
and therefore $2\leq U_\delta \le 3$ for
$ .061129446\ldots \leq \delta\leq .442695041\ldots$.
Using Maple VI we computed, for each
$u=2,2.1,\dots,2.9,3$, the value of
$\delta$ for which $U_\delta=u$, and then the
value of $I_\delta:=(1/U_\delta) \int_0^{U_\delta} \sigma_\delta(t) dt$:
\medskip

\centerline{\vbox{\offinterlineskip
\hrule
\halign{&\vrule \ # &\ \vrule \ # &\strut \ \vrule \ # &\ \vrule \ #\cr
\hfil $u=U_\delta$ \hfil &\hfil $\delta$ \hfil & \hfil $I_\delta$ \hfil &\cr
\noalign{\hrule}\cr
\noalign{\hrule}
\hfil$e^{1/2}$\hfil&\hfil1&\hfil.786938680\hfil &\cr
\hfil1.7\hfil&\hfil.884558536&\hfil.775994691\hfil&\cr
\hfil1.8\hfil&\hfil.701297528&\hfil.756132235\hfil&\cr
\hfil1.9\hfil&\hfil.557986983&\hfil.737993834\hfil&\cr
\hfil2.0\hfil&\hfil.442695041&\hfil.721347520\hfil&\cr
\hfil2.1\hfil&\hfil.353609191&\hfil.704809423\hfil&\cr
\hfil2.2\hfil&\hfil.286811221&\hfil.687757393\hfil&\cr
\hfil2.3\hfil&\hfil.234862762&\hfil.670734398\hfil&\cr
\hfil2.4\hfil&\hfil.193426306&\hfil.653994521\hfil&\cr
\hfil2.5\hfil&\hfil.159779207&\hfil.637653381\hfil&\cr
\hfil2.6\hfil&\hfil.132117433&\hfil.621755226\hfil&\cr
\hfil2.7\hfil&\hfil.109195664&\hfil.606305666\hfil&\cr
\hfil2.8\hfil&\hfil.090126952&\hfil.591288678\hfil&\cr
\hfil2.9\hfil&\hfil.074264622&\hfil.576675773\hfil&\cr
\hfil3.0\hfil&\hfil.061129446&\hfil.562431034\hfil&\cr}
\hrule}}

Using Maple, for $3\leq k\leq 17$, we give the lower bounds on
$\gamma_k$ (and $\gamma(S_k)$) that arise from
the above proof.

\centerline{\vbox{\offinterlineskip
\hrule
\halign{&\vrule \ # &\ \vrule \, \vrule\ # & \ \vrule \ #&\ \vrule \ # &\strut \ \vrule \ # &\ \vrule #\cr
\hfil $k$ \hfil &\hfil $\delta$ \hfil & \hfil $U_\delta$ \hfil & \hfil $\gamma_k\geq I_\delta=$ \hfil & \hfil $\gamma(S_k)\geq$ \hfil &\cr
\noalign{\hrule}\cr
\noalign{\hrule}
\hfil4\hfil&\hfil.3333333333&\hfil2.127612763&\hfil.7002748427&\hfil{.}\hfil&\cr
\hfil5\hfil&\hfil.2500000000&\hfil2.268355860&\hfil.6773393732&\hfil.7682091384\hfil&\cr
\hfil6\hfil&\hfil.2000000000&\hfil2.382637377&\hfil.6601481027&\hfil{.}\hfil&\cr
\hfil7\hfil&\hfil.1666666667&\hfil2.477839089&\hfil.6471915206&\hfil.7776179102\hfil&\cr
\hfil8\hfil&\hfil.1428571429&\hfil2.558879516&\hfil.6372773420&\hfil{.}\hfil&\cr
\hfil9\hfil&\hfil.1250000000&\hfil2.629113171&\hfil.6295761905&\hfil.7813572891\hfil&\cr
\hfil10\hfil&\hfil.1111111111&\hfil2.690898725&\hfil.6235174605&\hfil{.}\hfil&\cr
\hfil11\hfil&\hfil.1000000000&\hfil2.745943649&\hfil.6187030892&\hfil.7832215162\hfil&\cr
\hfil12\hfil&\hfil.09090909091&\hfil2.795516633&\hfil.6148498476&\hfil{.}\hfil&\cr
\hfil13\hfil&\hfil.08333333333&\hfil2.840582242&\hfil.6117521137&\hfil.7842851149\hfil&\cr
\hfil14\hfil&\hfil.07692307692&\hfil2.881888814&\hfil.6092577703&\hfil{.}\hfil&\cr
\hfil15\hfil&\hfil.07142857143&\hfil2.920027494&\hfil.6072523556&\hfil.7849492382\hfil&\cr
\hfil16\hfil&\hfil.06666666667&\hfil2.955472829&\hfil.6056484342&\hfil{.}\hfil&\cr
\hfil17\hfil&\hfil.06250000000&\hfil2.988611474&\hfil.6043783304&\hfil.7853917172\hfil&\cr}
\hrule}}

\enddocument